\documentclass{article}

\usepackage {amsmath, amsthm, amssymb}
\usepackage [top=1.5in, right=1in, bottom=1.5in, left=1in] {geometry}
\usepackage {tikz}
\usepackage {youngtab}
\usepackage{hyperref}
\hypersetup{
    colorlinks = true,
    citecolor = blue,
}

\usepackage {libertine}
\usepackage [T1]{fontenc}

\theoremstyle {Theorem}
\newtheorem {thm}  {Theorem}[subsection]
\newtheorem* {thm*}  {Theorem}
\newtheorem {prop} {Proposition}[subsection]

\newtheorem {cor}  {Corollary}[subsection]

\theoremstyle {definition}
\newtheorem {ex} {Example}[subsection]
\newtheorem {defn} {Definition}[subsection]
\newtheorem {rmk} {Remark}[subsection]

\title {$q$-Rationals and Finite Schubert Varieties}
\author {Nicholas Ovenhouse}
\date{}

\begin {document}

\maketitle

\begin {abstract}
    The classical $q$-analogue of the integers was recently generalized by Morier-Genoud and Ovsienko
    to give $q$-analogues of rational numbers. Some combinatorial interpretations are already known,
    namely as the rank generating functions for certain partially ordered sets. We review some of these
    interpretations, and additionally give a slightly novel approach in terms of planar graphs called snake graphs.
    Using the snake graph approach, we show that the numerators of $q$-rationals count the sizes of certain varieties over finite fields,
    which are unions of open Schubert cells in some Grassmannian.
\end {abstract}

\tableofcontents

\section {Background} \label{sec:background}

\bigskip

The classical ``\emph{$q$-analogue}'' of a positive integer $n$ is the following polynomial in $\Bbb{Z}[q]$
\[ [n]_q := 1 + q + q^2 + \cdots + q^{n-1} \]
It satisfies the property that evaluation at $q=1$ yields the number $n$. In general, any
``\emph{$q$-analogue}'' should satisfy at least this simple property --- that
evaluating $q=1$ (or more generally taking a limit $q \to 1$) should recover the classical case.

\bigskip

Some simple well-known examples are built from the $q$-integers in the obvious way. First is the $q$-factorial:
\[ [n]_q ! := [n]_q [n-1]_q \cdots [2]_q [1]_q \]
Also we have the $q$-binomial coefficients:
\[ \binom{n}{k}_q := \frac{[n]_q !}{[k]_q! [n-k]_q !} \]
However, the significance of a $q$-analogue would be superficial if the only connection were this ``$q=1$'' property.
Usually, a $q$-analogue satisfies some other interesting properties, and has some deeper connection with
the underlying expression. In particular, $q$-analogues often have two other interesting combinatorial interpretations.
\begin {itemize}
    \item[(1)] They count the size of some algebraic variety defined over a finite field $\Bbb{F}_q$.
    \item[(2)] They appear as weight generating functions of some combinatorial set for some natural statistic.
\end {itemize}

Recently, Morier-Genoud and Ovsienko defined $q$-analogues of rational numbers \cite{mgo}. Given a rational number $\frac{r}{s}$,
they define a rational function $\left[\frac{r}{s}\right]_q = \frac{\mathcal{R}(q)}{\mathcal{S}(q)}$, which is defined
in terms of the continued fraction expansion of $\frac{r}{s}$. The formal definition will be given later, but here are some examples:
\[ \left[\frac{5}{2}\right]_q = \frac{1+2q+q^2+q^3}{1+q}, \quad \quad  \left[\frac{10}{7}\right]_q = \frac{1+q+2q^2+3q^3+2q^4+q^5}{1+q+2q^2+2q^3+q^4} \]
Morier-Genoud and Ovsienko give some combinatorial interpretations and formulas for the polynomials $\mathcal{R}(q)$, but they
pose the open problem of finding more. They say the following in \cite{mgo}:

\bigskip

\begin {center}
\parbox{0.8\linewidth}{
``\emph{It is a challenging problem to find more different combinatorial and geometric interpretations of the polynomials
$\mathcal{R}$ and $\mathcal{S}$. A natural question is to connect the polynomials $\mathcal{R}$ and $\mathcal{S}$ with counting
of points in varieties defined over the finite fields $\Bbb{F}_q$. This property would be similar to that of the Gaussian 
$q$-binomial coefficients.''
}}
\end {center}

\bigskip

The aim of this paper is to address this question. We will review some known interpretations of the polynomials $\mathcal{R}(q)$,
as well as present some new ones. In particular, we answer the question of giving an interpretation in terms of
varieties over finite fields. The organization of the rest of the paper is as follows.

\bigskip

In the remainder of \textbf{Section \ref{sec:background}}, we discuss, for $[n]_q!$ and $\binom{n}{k}_q$, 
the two points mentioned above, how $q$-analogues often have interpretations in terms of finite algebraic varieties as well as
weight generating functions.

\medskip

In \textbf{Section \ref{sec:q_rational_def}}, we give the definition of $\left[\frac{r}{s}\right]_q$,
the $q$-analogue of the rational number $\frac{r}{s}$, and discuss some combinatorial interpretations in terms of weight generating
functions which have already appeared in the literature. In this section we define ``\emph{snake graphs}'', and
the main combinatorial interpretation of $q$-rationals for our purposes will be as weight generating functions
for certain paths on these graphs. In this section we also define two partitions $\lambda$ and $\mu$ (with $\mu < \lambda$)
associated with a rational number, which depend on the continued fraction expansion. These partitions play an
important part in the statement of the main theorem later.

\medskip

In \textbf{Section \ref{sec:Schubert}}, we begin by reviewing the decomposition of the Grassmannian $\mathrm{Gr}_k(n)$
into ``\emph{open Schubert cells}''. These cells are indexed by partitions which fit inside a $k \times (n-k)$ rectangle.
For each partition $\lambda$, we write $\Omega_\lambda^{\circ}$ for the open Schubert cell.
The main result (see \textbf{Theorem \ref{thm:main}}) is the following.

\bigskip

\begin {thm*} 
    Let $\frac{r}{s}$ be a rational number, and $\left[\frac{r}{s}\right]_q = \frac{\mathcal{R}(q)}{\mathcal{S}(q)}$. 
    Also let $\mu$ and $\lambda$ be the partitions defined in \textbf{Section \ref{sec:q_rational_def}} associated with $\frac{r}{s}$.
    Then 
    \[ q^{|\mu|} \mathcal{R}(q) = \left| \bigcup_{\mu \leq \nu \leq \lambda} \Omega_\nu^\circ \right| \]
\end {thm*}

\bigskip

\subsection {Geometry Over Finite Fields} \label{sec:geometry_finite}

\bigskip

Now let us recall how the above-mentioned examples (the $q$-factorial and $q$-binomial coefficients) can be interpreted as
the sizes of certain algebraic varieties over finite fields. The results in this section are well-known (see \cite{stanley} for example).

\medskip

Before stating the results, let us establish some notation. Let $V$ be a vector space (over a field $\Bbb{K}$), say of dimension $n$. 
The \emph{(complete) flag variety}, which we will denote by $\mathrm{Fl}(V)$, is the set of all complete flags in $V$. By a complete flag we mean
a chain of nested subspaces
\[ 0 \subset V_1 \subset V_2 \subset \cdots \subset V_{n-1} \subset V \tag{$\dim(V_i) = i$} \]
In the case that $V = \Bbb{K}^n$, we will just write $\mathrm{Fl}(n)$. 

\medskip

Recall also that the set of $k$-dimensional subspaces of a vector space is called a \emph{Grassmannian}, and denoted $\mathrm{Gr}_k(V)$,
or if $V = \Bbb{K}^n$, just $\mathrm{Gr}_k(n)$.

\bigskip

\begin {thm} \label{thm:finite_variety_sizes}
    If $\Bbb{K} = \Bbb{F}_q$ is a finite field, then
    \begin {itemize}
    \item[$(a)$] $\displaystyle \left| \mathrm{Fl}(n) \right| = [n]_q !$
    \item[$(b)$] $\displaystyle \left| \mathrm{Gr}_k(n) \right| = \binom{n}{k}_q$
    \end {itemize}
\end {thm}

\bigskip

\begin {ex} \label{ex:projective_space}
    The projective space $\Bbb{P}^{n-1}$ is the special case $\mathrm{Gr}_1(n)$, of lines in $\Bbb{K}^n$.
    In particular, we have that the $q$-integer $[n]_q = \binom{n}{1}_q$ is the size of $\Bbb{P}^{n-1}$ over $\Bbb{F}_q$.
    This is also easy to see using the usual definition of projective space: $\Bbb{P}^{n-1} = (\Bbb{K}^n \setminus \{0\}) / \sim$,
    where $\sim$ is the equivalence relation $p \sim \alpha \, p$ for any non-zero $\alpha \in \Bbb{K}$. In the case $\Bbb{K} = \Bbb{F}_q$,
    we have $\left| \Bbb{K} \setminus \{0\} \right| = q^n - 1$, and the number of non-zero $\alpha$'s is $q-1$.
    And indeed we have $[n]_q = \frac{q^n-1}{q-1}$.
\end {ex}

\bigskip

\subsection {Weight Generating Functions} \label{sec:wt_fns}

\bigskip

We will now recall the interpretations of $q$-factorials and $q$-binomial coefficients as weight generating functions. As in the previous section,
all results presented here are well-known (see \cite{stanley} for a standard reference).

\medskip

First we will discuss the $q$-factorial. The natural set that comes to mind when one thinks of the number $n!$ is the symmetric group $S_n$,
and indeed the $q$-analogue $[n]_q !$ is a weight generating function for a statistic on this set.
Given a permutation $\sigma \in S_n$, an \emph{inversion} of $\sigma$ is a pair $(i,j)$ such that $i < j$ and $\sigma(i) > \sigma(j)$.
Define $\mathrm{inv}(\sigma)$ to be the total number of inversions.

\bigskip

\begin {thm}
    $\displaystyle [n]_q ! = \sum_{\sigma \in S_n} q^{\mathrm{inv}(\sigma)}$
\end {thm}

\bigskip

Next we turn to the $q$-binomial coefficients. When one thinks of the numbers $\binom{n}{k}$, one naturally thinks of
the set $\binom{[n]}{k}$, the set of subsets of $[n] = \{1,2,\dots,n\}$ of size $k$. However, for the present purposes
it is more natural to consider a couple other sets which are in bijection with $\binom{[n]}{k}$. 

\medskip

Consider the set $P_{n,k}$ of lattice paths from $(0,0)$ to $(n-k,k)$ which only take unit steps right or up. All such paths
must take a total of $n$ steps, $k$ of which are ``up''. So there is a natural bijection $\binom{[n]}{k} \to P_{n,k}$. For a path
$p \in P_{n,k}$, let $|p|$ be the area above $p$ and below the horizontal line $y=k$.

\bigskip

\begin {thm} \label{thm:q_binom_gen_fn}
    $\displaystyle \sum_{p \in P_{n,k}} q^{|p|} = \binom{n}{k}_q$
\end {thm}

\bigskip

This can also be described nicely in terms of integer partitions. Recall that a \emph{partition} $\lambda$ is a sequence
of weakly decreasing integers $\lambda_1 \geq \lambda_2 \geq \cdots \geq \lambda_k$. If $\sum \lambda_i = n$, we write $\lambda \vdash n$,
and also $|\lambda| = n$. Partitions are visualized by their Young diagram, which is an array of boxes (left-aligned) with $\lambda_i$
boxes in row $i$. Let $Y_{n,k}$ be the set of all partitions whose Young diagram fits inside the rectangle with height $k$ and width $n-k$.
There is a clear bijection $Y_{n,k} \to P_{n,k}$. The bottom boundary of any Young diagram is a lattice path in $P_{n,k}$, and the number
of boxes in the Young diagram is precicely the area above the lattice path.

\bigskip

\begin {cor} \label{cor:q_binom_gen_fn}
    $\displaystyle \sum_{\lambda \in Y_{n,k}} q^{|\lambda|} = \binom{n}{k}_q$
\end {cor}

\bigskip

Lastly, let us just mention that the previous expressions can also be thought of in terms of posets. In particular, the set of all Young diagrams (or partitions)
has a natural partial order where $\mu \leq \lambda$ means that $\mu_i \leq \lambda_i$ for all $i$, or equivalently the Young diagram of $\mu$ fits inside $\lambda$.
This poset is called \emph{Young's lattice}. The poset has a natural \emph{rank} function, given by the size of the partition (the number of boxes).
If we let $\lambda = (n-k)^k$ denote the Young diagram given by the $k \times (n-k)$ rectangle, then the set $Y_{n,k}$ above is simply the interval
of Young's lattice $[\varnothing, \lambda] = \{ \mu ~|~ \mu \leq \lambda\}$, and the generating function given above is the rank generating function for
this poset. This simply means the coefficient of $q^k$ is the number of elements of rank $k$. 

\bigskip

\begin {ex}
    The poset structure for $Y_{4,2}$ and $P_{4,2}$ is pictured in \textbf{Figure \ref{fig:young_poset}}. 
    The sizes of the ranks are the coefficients of $\binom{4}{2}_q = 1+q+2q^2+q^3+q^4$.
    \begin {figure}[h]
    \centering
    \begin {tikzpicture}[scale=0.3]
        \begin {scope}[shift={(0,4)}]
            \fill[gray] (0,1) -- (1,1) -- (1,2) -- (0,2) -- cycle;
        \end {scope}

        \begin {scope}[shift={(-3,7)}]
            \fill[gray] (0,0) -- (1,0) -- (1,2) -- (0,2) -- cycle;
        \end {scope}

        \begin {scope}[shift={(3,7)}]
            \fill[gray] (0,1) -- (2,1) -- (2,2) -- (0,2) -- cycle;
        \end {scope}

        \begin {scope}[shift={(0,10)}]
            \fill[gray] (0,0) -- (1,0) -- (1,1) -- (2,1) -- (2,2) -- (0,2) -- cycle;
        \end {scope}

        \begin {scope}[shift={(0,14)}]
            \fill[gray] (0,0) -- (2,0) -- (2,2) -- (0,2) -- cycle;
        \end {scope}

        \foreach \x/\y in {0/0, 0/4, -3/7, 3/7, 0/10, 0/14} {
            \draw (\x,\y) -- (\x+2,\y) -- (\x+2,\y+2) -- (\x,\y+2) -- cycle;
            \draw (\x+1,\y) -- (\x+1,\y+2);
            \draw (\x,\y+1) -- (\x+2,\y+1);
        }

        \begin {scope}[shift={(0,0)}]
            \draw[red, line width=1.2] (0,0) -- (0,2) -- (2,2);
        \end {scope}

        \begin {scope}[shift={(0,4)}]
            \draw[red, line width=1.2] (0,0) -- (0,1) -- (1,1) -- (1,2) -- (2,2);
        \end {scope}

        \begin {scope}[shift={(-3,7)}]
            \draw[red, line width=1.2] (0,0) -- (1,0) -- (1,2) -- (2,2);
        \end {scope}

        \begin {scope}[shift={(3,7)}]
            \draw[red, line width=1.2] (0,0) -- (0,1) -- (2,1) -- (2,2);
        \end {scope}

        \begin {scope}[shift={(0,10)}]
            \draw[red, line width=1.2] (0,0) -- (1,0) -- (1,1) -- (2,1) -- (2,2);
        \end {scope}

        \begin {scope}[shift={(0,14)}]
            \draw[red, line width=1.2] (0,0) -- (2,0) -- (2,2);
        \end {scope}
    \end {tikzpicture}
    \caption {Poset structure for $Y_{4,2}$}
    \label {fig:young_poset}
    \end {figure}
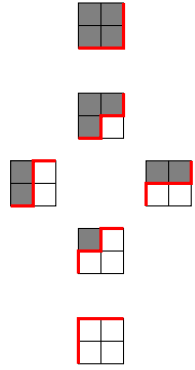
\end {ex}

\bigskip

\section {$q$-Rational Numbers} \label{sec:q_rational_def}

\bigskip

In \cite{mgo}, Morier-Genoud and Ovsienko extend the definition of $[n]_q$ to include the case when $n = \frac{r}{s} \in \Bbb{Q}$
is a rational number. Their definition uses the continued fraction expansion of $\frac{r}{s}$.
Specifically, if $\frac{r}{s} = [a_1,a_2,\dots,a_{2m}]$ is the finite continued fraction
expansion (and $a_1 > 1$), then they define
\[ 
    \left[ \frac{r}{s} \right]_q := 
    [a_1]_q + \cfrac{q^{a_1}}{
        [a_2]_{q^{-1}} + \cfrac{q^{-a_2}}{
            [a_3]_q + \cfrac{q^{a_3}}{
                [a_4]_{q^{-1}} + \cfrac{q^{-a_4}}{
                    \ddots \, + \, \cfrac{q^{a_{2m-1}}}{
                        [a_{2m}]_{q^{-1}}
                    } 
                }
            }
        }
    } 
\]
Morier-Genoud and Ovsienko also give several other ways to compute these expressions, including some determinantal formulas,
a recursive procedure involving triangulated polygons, and using products of $q$-deformed matrices in $\mathrm{PSL}_2(\Bbb{Z})$.
We will briefly describe this last interpretation in terms of matrices (for more details, see \cite{mgo} and \cite{mgl}). 

\medskip

Define the following two matrices:
\[ A = \begin{pmatrix} 1 & 1 \\ 0 & 1 \end{pmatrix} \quad \text{ and } \quad B = \begin{pmatrix} 1 & 0 \\ 1 & 1 \end{pmatrix} \]
They generate $\mathrm{PSL}_2(\Bbb{Z})$, and the group acts transitively on $\Bbb{Q} \cup \{\infty\}$ by the rule:
\[ \begin{pmatrix} a&b\\c&d \end{pmatrix} \cdot x := \frac{ax+b}{cx+d} \]
Define $q$-deformed versions of these matrices:
\[ A_q = \begin{pmatrix} q & 1 \\ 0 & 1 \end{pmatrix} \quad \text{ and } \quad B_q = \begin{pmatrix} q & 0 \\ q & 1 \end{pmatrix} \]
Let $\mathrm{PSL}_2^q(\Bbb{Z})$ be the group generated by $A_q$ and $B_q$, modulo scaling by monomials $q^{\pm n}$. In general, for
any $M \in \mathrm{PSL}_2(\Bbb{Z})$, we get a corresponding element $[M]_q \in \mathrm{PSL}_2^q(\Bbb{Z})$. This group
acts on the set of power series $\Bbb{Z}[[q]] \cup \{\infty\}$ by the same formula as above.
The following is another, more conceptual, definition of $q$-rational numbers.

\bigskip

\begin {thm} \label{thm:q_matrices} \cite{mgo} \cite{mgl} The $q$-deformation $x \mapsto [x]_q$ commutes with the $\mathrm{PSL}_2(\Bbb{Z})$ action.
    In particular, this means if $x \in \Bbb{Q}$ and $M \in \mathrm{PSL}_2(\Bbb{Z})$, then $[M \cdot x]_q = [M]_q \cdot [x]_q$.
\end {thm}

\bigskip

\begin {ex} \label{ex:seven_thirds_computed}
    The continued fraction expansion of $\frac{7}{3}$ is $[2,3]$. Using the definition, we get
    \[ \left[\frac{7}{3}\right]_q = [2]_q + \cfrac{q^2}{[3]_{q^{-1}}} = (1+q) + \frac{q^2}{1+q^{-1}+q^{-2}} = \frac{1+2q+2q^2+q^3+q^4}{1+q+q^2} \] 
    Alternatively, we can use \textbf{Theorem \ref{thm:q_matrices}} to compute. Notice that the generators $A$ and $B$
    act by $A \cdot x = x+1$ and $B \cdot x = \frac{x}{x+1}$. Notice that $\frac{7}{3} = A^2B^2 \cdot 1$:
    \[ 1 \xrightarrow{~B~} \frac{1}{2} \xrightarrow{~B~} \frac{1}{3} \xrightarrow{~A~} \frac{4}{3} \xrightarrow{~A~} \frac{7}{3} \]
    According to \textbf{Theorem \ref{thm:q_matrices}}, we should be able to start with $1$ and apply $A_q^2B_q^2$ to compute $\left[ \frac{7}{3}\right]_q$.
    Notice that the action of $A_q$ and $B_q$ are given by $A_q \cdot x = 1 + qx$ and $B_q \cdot x = \frac{qx}{1+qx}$.
    \[ 1 \xrightarrow{~B_q~} \frac{q}{1+q} \xrightarrow{~B_q~} \frac{q^2}{1+q+q^2} \xrightarrow{~A_q~} \frac{1+q+q^2+q^3}{1+q+q^2} \xrightarrow{~A_q~} \frac{1+2q+2q^2+q^3+q^4}{1+q+q^2} \]
\end {ex}

\bigskip

\subsection {Combinatorial Interpretations} \label{sec:q_rational_wt_gen}

\bigskip

Let $\left[ \frac{r}{s} \right]_q = \frac{\mathcal{R}(q)}{\mathcal{S}(q)}$. One of the interpretations given in \cite{mgo} says that
the coefficient of $q^k$ in $\mathcal{R}(q)$ is the number of ``\emph{$k$-vertex closures}'' in a certain directed graph. Although
worded differently, this is equivalent to the following interpretation in terms of posets (see \cite{mss}). 

\bigskip

\begin {defn} 
    Define a poset $F\left(\frac{r}{s}\right)$ on the set $\{x_1,x_2,\dots,x_{N-1}\}$, where $N = \sum_{i=1}^{2m} a_i$, with cover relations:
    \[ x_1 < x_2 < \cdots < x_{a_1} > x_{a_1+1} > \cdots > x_{a_1+a_2} < \cdots \] 
    In other words, the Hasse diagram is a ``\emph{fence}'' which goes up $(a_1-1)$ times,
    then down $a_2$ times, then up $a_3$ times, $\dots$, and finally up $a_{2m-1}$ times, and down $(a_{2m}-1)$ times. 
\end {defn}

\bigskip

Let $L\left(\frac{r}{s}\right)$ be the lattice of order ideals in $F\left(\frac{r}{s}\right)$.
In this language, what Morier-Genoud and Ovsienko call a ``\emph{$k$-vertex closure}'' is the same as a lower order ideal of $F\left(\frac{r}{s}\right)$.
In particular we have the following. 

\bigskip

\begin {thm} \cite{mgo} 
    $\mathcal{R}(q)$ is the rank generating function of $L\left(\frac{r}{s}\right)$.
\end {thm}

\bigskip

\begin {rmk}
    The connection between these posets and $F$-polynomials in cluster algebras has been noted
    in \cite{msw_11}, \cite{rabideau}, \cite{bg_21}, and \cite{claussen}. These posets were called ``\emph{fence posets}'' in \cite{mss}
    and ``\emph{piece-wise linear posets}'' in \cite{bg_21}.
    It was shown in \cite{mgo} that $\mathcal{R}(q)$ can be obtained from an $F$-polynomial by specialization of the variables.
\end {rmk}
    
\bigskip

\begin {ex} \label{ex:seven_thirds}
    We continue with the running example of $\frac{r}{s} = \frac{7}{3}$ from \textbf{Example \ref{ex:seven_thirds_computed}}.
    Recall that the continued fraction for $\frac{7}{3}$ is $[2,3]$. The corresponding fence $F(7/3)$ is
    
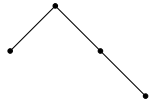
\begin {center}
    \begin {tikzpicture}[scale=0.6]
        \draw[fill=black] (0,0) circle (0.05);
        \draw[fill=black] (1,1) circle (0.05);
        \draw[fill=black] (2,0) circle (0.05);
        \draw[fill=black] (3,-1) circle (0.05);

        \draw (0,0) -- (1,1) -- (3,-1);
    \end {tikzpicture}
    \end {center}
    The Lattice $L(7/3)$ of order ideals in $F(7/3)$ is pictured in \textbf{Figure \ref{fig:seven_thirds}}. 
    Notice the number of elements of rank $k$ is the coefficient of $q^k$ in $\mathcal{R}(q) = 1+2q+2q^2+q^3+q^4$.
    \begin {figure}[h]
    \centering
    \begin {tikzpicture} [scale=0.25]
        \foreach \x/\y in {0/0, -5/4, 5/4, 0/8, 10/8, 5/12, 5/16} {
            \draw[fill=black] (\x,\y) circle (0.05);
            \draw[fill=black] (\x+1,\y+1) circle (0.05);
            \draw[fill=black] (\x+2,\y+0) circle (0.05);
            \draw[fill=black] (\x+3,\y-1) circle (0.05);
            \draw (\x,\y) -- (\x+1,\y+1) -- (\x+3,\y-1);
        }

        \begin {scope} [shift = {(0,0)}]
        \end {scope} 

        \begin {scope} [shift = {(-5,4)}]
            \draw[red] (0,0) circle (0.3);
        \end {scope} 

        \begin {scope} [shift = {(5,4)}]
            \draw[red] (3,-1) circle (0.3);
        \end {scope} 

        \begin {scope} [shift = {(0,8)}]
            \draw[red] (0,0) circle (0.3);
            \draw[red] (3,-1) circle (0.3);
        \end {scope} 

        \begin {scope} [shift = {(10,8)}]
            \draw[red] (3,-1) circle (0.3);
            \draw[red] (2,0) circle (0.3);
        \end {scope} 

        \begin {scope} [shift = {(5,12)}]
            \draw[red] (3,-1) circle (0.3);
            \draw[red] (2,0) circle (0.3);
            \draw[red] (0,0) circle (0.3);
        \end {scope} 

        \begin {scope} [shift = {(5,16)}]
            \draw[red] (3,-1) circle (0.3);
            \draw[red] (2,0) circle (0.3);
            \draw[red] (0,0) circle (0.3);
            \draw[red] (1,1) circle (0.3);
        \end {scope} 

        \draw (-0.5,1.5) -- (-1.5,2.5); 
        \draw (4.5,5.5) -- (3.5,6.5);
        \draw (9.5,9.5) -- (8.5,10.5);
        \draw (3.5,1.5) -- (4.5,2.5);
        \draw (8.5,5.5) -- (9.5,6.5);
        \draw (-1.5,5.5) -- (-0.5,6.5);
        \draw (3.5,9.5) -- (4.5,10.5);
        \draw (6.5,13.5) -- (6.5,15);

    \end {tikzpicture}
    \caption {Poset structure of $L(7/3)$}
    \label {fig:seven_thirds}
    \end {figure}
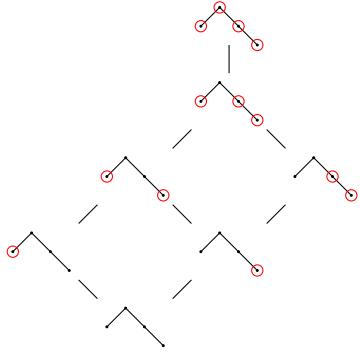
\end {ex}

\bigskip

\subsection {Snake Graphs}

\bigskip

Yet another equivalent description can be given in terms of ``\emph{snake graphs}''. These are certain types of planar graphs built out
of square tiles, such that each tile is either above or to the right of the previous one. To each snake graph $G$ we can naturally
associate a \emph{word} $W(G)$ in the alphabet $\{R,U\}$, indicating whether going from one tile to the next is ``\emph{right}'' or ``\emph{up}''.
Some examples are shown in \textbf{Figure \ref{fig:snakes}}.

\begin {figure}[h]
\centering
\begin {tikzpicture}[scale=0.8]
    \draw (0,0) -- ++(3,0) -- ++(0,1) -- ++(-3,0) -- cycle;
    \draw (1,0) -- (1,1);
    \draw (2,0) -- (2,1);

    \draw (1.5,-2) node {$W(G) = RR$};
    \draw (1.5,-3) node {$G = G(4/1)$};
    \draw (1.5,-4) node {$4/1 = [3,1]$};

    \begin {scope}[shift={(5,-1)}]
        \draw (0,0) -- (1,0) -- (1,1) -- (2,1) -- (2,2) -- (3,2) -- (3,3) -- (1,3) -- (1,2) -- (0,2) -- cycle;
        \draw (0,1) -- (1,1) -- (1,2) -- (2,2) -- (2,3);

        \draw (1.5,-1) node {$W(G) = URUR$};
        \draw (1.5,-2) node {$G = G(13/8)$};
        \draw (1.5,-3) node {$13/8 = [1,1,1,1,1,1]$};
    \end {scope}

    \begin {scope}[shift={(10,-1)}]
        \draw (0,0) -- (2,0) -- (2,1) -- (4,1) -- (4,2) -- (1,2) -- (1,1) -- (0,1) -- cycle;
        \draw (1,0) -- (1,1) -- (2,1) -- (2,2);
        \draw (3,1) -- (3,2);

        \draw (1.5,-1) node {$W(G) = RURR$};
        \draw (1.5,-2) node {$G = G(11/4)$};
        \draw (1.5,-3) node {$11/4 = [2,1,2,1]$};
    \end {scope}
\end {tikzpicture}
\caption {Examples of snake graphs}
\label {fig:snakes}
\end {figure}

\medskip

\begin {defn}
    Let $\overline{R} = U$ and $\overline{U} = R$. If $w = w_1 w_2 \cdots w_k$ is a word in the alphabet $\{R,U\}$,
    define the \emph{dual word} to be $w^* = \overline{w}_1 w_2 \overline{w}_3 w_4 \cdots \overline{w}_{2n-1} w_{2n}$ or 
    $w^* = \overline{w}_1 w_2 \overline{w}_3 w_4 \cdots w_{2n} \overline{w}_{2n+1}$ (depending on parity). In other words, $w^*$ is obtained
    by toggling the odd-indexed letters. Extend this definition to snake graphs by defining $G^*$ to be the snake graph so that $W(G^*) = W(G)^*$.
\end {defn}

\bigskip

\begin {rmk}
    Snake graphs have appeared numerous times
    in the cluster algebra literature (see for example \cite{cs_18} \cite{rabideau} \cite{msw_11}, \cite{propp}). 
    The dual construction ($G^*$ rather than $G$) seems to be more popular in the literature 
    (in \cite{propp}, Propp first constructs $G^*$ and calls $G$ the ``dual'' snake). In what follows, results from the literature
    are paraphrased in terms of $G$ rather than $G^*$.
\end {rmk}

\bigskip

\begin {defn} \cite{cs_18}
    Given a continued fraction $\frac{r}{s} = [a_1,\dots,a_{2m}]$, define a snake graph $G\left( \frac{r}{s} \right)$ such that
    its word is $W(G) = R^{a_1-1}U^{a_2}R^{a_3}U^{a_4} \cdots R^{a_{2m-1}} U^{a_{2m}-1}$. 
    \textbf{Figure \ref{fig:snakes}} also indicates which rational number $\frac{r}{s}$ is associated to each snake graph.
\end {defn}

\bigskip

\bigskip

\begin {defn}
    Given a snake graph $G$, let $P(G)$ be the set of lattice paths on $G$,
    going from the bottom-left corner to the top-right corner, which only take steps right or up.
    We also use the notation $P(r/s) = P(G(r/s))$.
\end {defn}

\bigskip
   
\begin {thm} \label{thm:path_count} \cite{propp} \cite{cs_18} 
    The number of lattice paths on $G(r/s)$ is $r$. In other words, $|P(r/s)| = r$.
\end {thm}

\bigskip 

We now give a $q$-analogue of \textbf{Theorem \ref{thm:path_count}}, where $r$
is replaced by $\mathcal{R}(q)$, the numerator of the corresponding $q$-rational.
Note that there is a natural partial order on $P(G)$ so that $a \leq b$ if the set of boxes
above $a$ is a subset of the boxes above $b$. Let $|a|$ denote the number of boxes
above the path $a$.

\bigskip

\begin {figure}[h]
\centering
\begin {tikzpicture}[scale=0.4]
    \draw[fill=black] (0,0) circle (0.05);
    \draw[fill=black] (1,1) circle (0.05);
    \draw[fill=black] (2,0) circle (0.05);
    \draw[fill=black] (3,1) circle (0.05);
    \draw[fill=black] (4,2) circle (0.05);
    \draw[fill=black] (5,3) circle (0.05);
    \draw[fill=black] (6,2) circle (0.05);
    \draw[fill=black] (7,1) circle (0.05);

    \draw (0,0) -- (1,1) -- (2,0) -- (5,3) -- (7,1);

    \draw (3.5,-2) node {$F(r/s)$};

    \draw[->] (8,1.5) -- (9,1.5);

    \begin {scope}[shift={(10,3)}]
        \draw[fill=black] (0,0) circle (0.05);
        \draw[fill=black] (1,-1) circle (0.05);
        \draw[fill=black] (2,0) circle (0.05);
        \draw[fill=black] (3,-1) circle (0.05);
        \draw[fill=black] (4,-2) circle (0.05);
        \draw[fill=black] (5,-3) circle (0.05);
        \draw[fill=black] (6,-2) circle (0.05);
        \draw[fill=black] (7,-1) circle (0.05);

        \draw (0,0) -- (1,-1) -- (2,0) -- (5,-3) -- (7,-1);

        \draw (3.5,-5) node {reflect};

        \draw[->] (8,-1.5) -- (9,-1.5);
    \end {scope}

    \begin {scope}[shift={(20,0)}]
        \draw[fill=black] (0,0) circle (0.05);
        \draw[fill=black] (1,0) circle (0.05);
        \draw[fill=black] (1,1) circle (0.05);
        \draw[fill=black] (2,1) circle (0.05);
        \draw[fill=black] (3,1) circle (0.05);
        \draw[fill=black] (4,1) circle (0.05);
        \draw[fill=black] (4,2) circle (0.05);
        \draw[fill=black] (4,3) circle (0.05);

        \draw (0,0) -- (1,0) -- (1,1) -- (4,1) -- (4,3);

        \draw (2,-2) node {rotate};

        \draw[->] (6,1.5) -- (7,1.5);
    \end {scope}

    \begin {scope}[shift={(30,0)}]
        \draw[gray,fill=gray] (0,0) circle (0.05);
        \draw[gray,fill=gray] (1,0) circle (0.05);
        \draw[gray,fill=gray] (1,1) circle (0.05);
        \draw[gray,fill=gray] (2,1) circle (0.05);
        \draw[gray,fill=gray] (3,1) circle (0.05);
        \draw[gray,fill=gray] (4,1) circle (0.05);
        \draw[gray,fill=gray] (4,2) circle (0.05);
        \draw[gray,fill=gray] (4,3) circle (0.05);

        \draw[gray] (0,0) -- (1,0) -- (1,1) -- (4,1) -- (4,3);

        \draw (2,-2) node {$G(r/s)$};

        \foreach \x\y in {0/0, 1/0, 1/1, 2/1, 3/1, 4/1, 4/2, 4/3} {
            \draw (\x-0.5,\y-0.5) -- (\x+0.5,\y-0.5) -- (\x+0.5,\y+0.5) -- (\x-0.5,\y+0.5) -- cycle;
        }
    \end {scope}
\end {tikzpicture}
\caption {Proof of \textbf{Proposition \ref{prop:rank_function}}}
\label {fig:snake_core}
\end {figure}
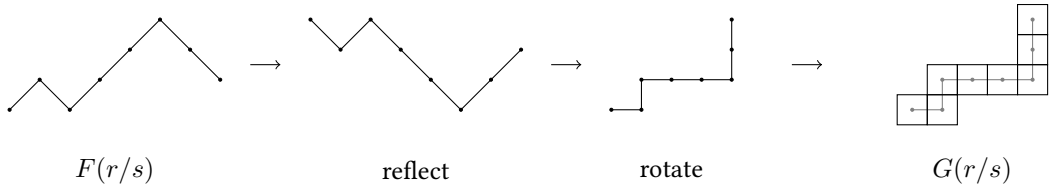

\begin {prop} \label{prop:rank_function}
    $\mathcal{R}(q)$ is the rank-generating function of the poset $P(r/s)$. In other words,
    \[ \sum_{p \in P(r/s)} q^{|p|} = \mathcal{R}(q) \]
\end {prop}
\begin {proof}
    It is easy to see that $P(r/s)$ is isomorphic to the poset $L(r/s)$ described earlier.
    The points of $F(r/s)$ correspond to the boxes of the snake graph, and
    the order ideals of $F(r/s)$ are the sets of boxes above the lattice paths.

    More specifically, we can construct $G(r/s)$ from $F(r/s)$ as follows (see \textbf{Figure \ref{fig:snake_core}} for an illustration). 
    Start with the Hasse diagram of $F(r/s)$, and reflect it over a horizontal line. Next, rotate it $45^\circ$ counter-clockwise.
    Then draw a square around each vertex. This will be $G(r/s)$.
\end {proof}

\bigskip

\begin {ex} \label{ex:seven_thirds_again}
    The poset $P(7/3)$ is shown in \textbf{Figure \ref{fig:seven_thirds_again}}. 
    Compare with $L(7/3)$ from \textbf{Example \ref{ex:seven_thirds}} and \textbf{Figure \ref{fig:seven_thirds}}.

    \begin {figure}[h]
    \centering
    \begin {tikzpicture} [scale=0.25]

        \begin {scope} [shift = {(-5,4)}]
            \fill[gray] (0,0) -- (1,0) -- (1,1) -- (0,1) -- cycle;
        \end {scope} 

        \begin {scope} [shift = {(5,4)}]
            \fill[gray] (1,2) -- (2,2) -- (2,3) -- (1,3) -- cycle;
        \end {scope} 

        \begin {scope} [shift = {(0,8)}]
            \fill[gray] (0,0) -- (1,0) -- (1,1) -- (0,1) -- cycle;
            \fill[gray] (1,2) -- (2,2) -- (2,3) -- (1,3) -- cycle;
        \end {scope} 

        \begin {scope} [shift = {(10,8)}]
            \fill[gray] (1,1) -- (2,1) -- (2,3) -- (1,3) -- cycle;
        \end {scope} 

        \begin {scope} [shift = {(5,12)}]
            \fill[gray] (1,1) -- (2,1) -- (2,3) -- (1,3) -- cycle;
            \fill[gray] (0,0) -- (1,0) -- (1,1) -- (0,1) -- cycle;
        \end {scope} 

        \begin {scope} [shift = {(5,18)}]
            \fill[gray] (0,0) -- (2,0) -- (2,3) -- (1,3) -- (1,1) -- (0,1) -- cycle;
        \end {scope} 

        \foreach \x/\y in {0/0, -5/4, 5/4, 0/8, 10/8, 5/12, 5/18} {
            \draw (\x,\y) -- (\x+2,\y) -- (\x+2,\y+3) -- (\x+1,\y+3) -- (\x+1,\y+1) -- (\x,\y+1) -- cycle;
            \draw (\x+1,\y) -- (\x+1,\y+1) -- (\x+2,\y+1);
            \draw (\x+1,\y+2) -- (\x+2,\y+2);
        }

        \begin {scope} [shift = {(0,0)}]
            \draw[red, line width=1] (0,0) -- (0,1) -- (1,1) -- (1,3) -- (2,3);
        \end {scope} 

        \begin {scope} [shift = {(-5,4)}]
            \draw[red, line width=1] (0,0) -- (1,0) -- (1,3) -- (2,3);
        \end {scope} 

        \begin {scope} [shift = {(5,4)}]
            \draw[red, line width=1] (0,0) -- (0,1) -- (1,1) -- (1,2) -- (2,2) -- (2,3);
        \end {scope} 

        \begin {scope} [shift = {(0,8)}]
            \draw[red, line width=1] (0,0) -- (1,0) -- (1,1) -- (1,2) -- (2,2) -- (2,3);
        \end {scope} 

        \begin {scope} [shift = {(10,8)}]
            \draw[red, line width=1] (0,0) -- (0,1) -- (2,1) -- (2,3);
        \end {scope} 

        \begin {scope} [shift = {(5,12)}]
            \draw[red, line width=1] (0,0) -- (1,0) -- (1,1) -- (2,1) -- (2,3);
        \end {scope} 

        \begin {scope} [shift = {(5,18)}]
            \draw[red, line width=1] (0,0) -- (2,0) -- (2,3);
        \end {scope} 

        \draw (-0.5,1.5) -- (-1.5,2.5); 
        \draw (4.5,5.5) -- (3.5,6.5);
        \draw (9.5,9.5) -- (8.5,10.5);
        \draw (3.5,1.5) -- (4.5,2.5);
        \draw (8.5,5.5) -- (9.5,6.5);
        \draw (-1.5,5.5) -- (-0.5,6.5);
        \draw (3.5,9.5) -- (4.5,10.5);
        \draw (6.5,16) -- (6.5,17.5);

    \end {tikzpicture}
    \caption {Poset structure of $P(7/3)$}
    \label {fig:seven_thirds_again}
    \end {figure}
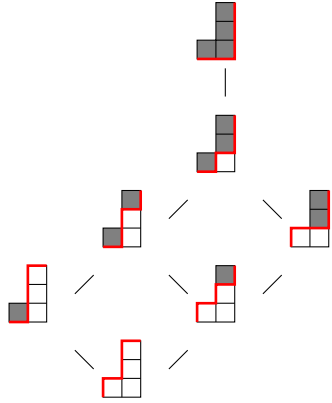
\end {ex}

\bigskip

The snake graph interpretation allows to make a nice analogy with the situation for $q$-binomial coefficients. In
particular, $P(r/s)$ is also isomorphic to an interval in Young's lattice. This follows from the observation that
any snake graph is a skew Young diagram of shape $\lambda/\mu$ where $\lambda$ is determined by the bottom boundary
of $G$, and $\mu$ is determined by the top boundary of $G$.

\bigskip

\begin {defn}
    If $G$ is a snake graph, let $\lambda(G)$ be the partition such that $G$ is the subset of boxes adjacent to the
    bottom boundary of $\lambda$. Also, let $\mu(G) < \lambda(G)$ be the partition such that $G$ is the skew
    Young diagram of shape $\lambda(G)/\mu(G)$. If $G = G(r/s)$, then we use the notations $\lambda(r/s)$ and $\mu(r/s)$.
\end {defn}

\bigskip

\begin {ex}
    Some examples of $\lambda$ and $\mu$, along with their corresponding fractions $\frac{r}{s}$, are shown in
    \textbf{Figure \ref{fig:lambda_mu}}.
    \begin {figure}[hb]
    \centering
    \begin {tikzpicture}[scale=0.6]
        \draw (0,0) -- ++(3,0) -- ++(0,1) -- ++(-3,0) -- cycle;
        \draw (1,0) -- (1,1);
        \draw (2,0) -- (2,1);

        \draw (1.5,-2) node {$\frac{r}{s} = \frac{4}{1} = [3,1]$};
        \draw (1.5,-3) node {$W(4/1) = RR$};
        \draw (1.5,-4) node {$\lambda = (3)$, $\mu=\varnothing$};

        \begin {scope}[shift={(8,-1)}]
            \fill[gray, opacity=0.5] (0,1) -- (1,1) -- (1,3) -- (0,3) -- cycle;
            \draw (0,0) -- (2,0) -- (2,2) -- (3,2) -- (3,3) -- (0,3) -- cycle;
            \draw (1,0) -- (1,3);
            \draw (2,2) -- (2,3);
            \draw (0,1) -- (2,1);
            \draw (0,2) -- (2,2);

            \draw (1.5,-1) node {$\frac{r}{s} = \frac{12}{5} = [2,2,1,1]$};
            \draw (1.5,-2) node {$W(12/5) = RUUR$};
            \draw (1.5,-3) node {$\lambda=(3,2^2)$, $\mu=(1^2)$};
        \end {scope}

        \begin {scope}[shift={(16,-1)}]
            \fill[gray, opacity=0.5] (0,2) -- (2,2) -- (2,3) -- (3,3) -- (3,4) -- (0,4) -- cycle;
            \draw (0,0) -- (1,0) -- (1,1) -- (3,1) -- (3,2) -- (4,2) -- (4,4) -- (0,4) -- cycle;
            \draw (0,1) -- (1,1) -- (1,4);
            \draw (0,2) -- (3,2) -- (3,4);
            \draw (0,3) -- (4,3);
            \draw (2,1) -- (2,4);

            \draw (1.5,-1) node {$\frac{r}{s} = \frac{31}{18} = [1,1,2,1,1,2]$};
            \draw (1.5,-2) node {$W(31/18) = URRURU$};
            \draw (1.5,-3) node {$\lambda=(4^2,3,1)$, $\mu=(3,2)$};
        \end {scope}
    \end {tikzpicture}
    \caption {Examples of $\lambda(r/s)$ and $\mu(r/s)$. The partition $\mu$ is shaded in gray.}
    \label {fig:lambda_mu}
    \end {figure}
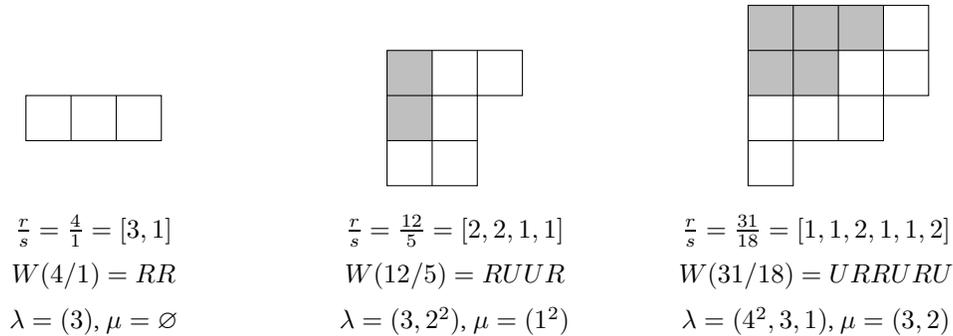
\end {ex}

\bigskip

In order to state the next result, we introduce the following notation. We will write partitions as
$\lambda = (\lambda_1^{b_1},\lambda_2^{b_2}, \dots)$, where $\lambda_i^{b_i}$ means that $\lambda_i$ is
repeated $b_i$ times. For example, $(3^4,2,1^2) = (3,3,3,3,2,1,1)$.

\bigskip

\begin {prop}
    Suppose $\frac{r}{s} = [a_1,b_1,\dots,a_m,b_m]$. The partitions $\lambda(r/s)$ and $\mu(r/s)$ are given explicitly as follows.
    \begin {itemize}
        \item[$(a)$] $\lambda = (\lambda_1^{b_m},\lambda_2^{b_{m-1}},\dots,\lambda_m^{b_1})$, where $\lambda_k = \sum\limits_{i=1}^{m+1-k} a_i$.
        \item[$(b)$] $\mu = (\mu_1^{c_1}, \dots, \mu_m^{c_m})$, where $\mu_i = \lambda_i-1$ and $c_i$ are given by $c_1=b_m-1$ and $c_i=b_{m+1-i}$ otherwise.
    \end {itemize}
\end {prop}
\begin {proof}
    $(a)$ We will induct on $m$. Recall the definition of the snake graph $G(r/s)$. It is defined by the word $W(r/s) = R^{a_1-1}U^{b_1}R^{a_2}U^{b_2} \cdots R^{a_m}U^{b_m-1}$.
    By definition of $\lambda(r/s)$, the partition $\lambda$ will have the same bottom boundary as $G(r/s)$. Suppose $m=1$.
    Then $W(r/s) = R^{a_1-1}U^{b_1-1}$. The young diagram which has this as its bottom boundary is the rectangle with width $a_1$ and height $b_1$ (i.e. $\lambda = (a_1^{b_1})$).

    \medskip

    Now suppose the result is true for all continued fractions $\frac{r'}{s'} = [a_1,b_1,\dots,a_{m-1},b_{m-1}]$.
    Then by induction $\lambda(r'/s') = \lambda' = (\lambda_1^{b_{m-1}},\lambda_2^{b_{m-2}}, \dots, \lambda_{m-1}^{b_1})$, where
    $\lambda_k = \sum_{i=1}^{m-k} a_i$. The snake graph for $\frac{r}{s} = [a_1,b_1,\dots,a_m,b_m]$ is build from the previous one
    by adding some more boxes to the end. This means the parition $\lambda = \lambda(r/s)$ will be obtained from $\lambda' = \lambda(r'/s')$
    by adding some more rows to the top.

    \medskip

    First let us point out a subtle point. The word $W(r'/s')$ ends with $U^{b_{m-1}-1}$, but $W(r/s)$ contains $U^{b_{m-1}}$.
    This means we must add one more box going up at the end of $G(r'/s')$ before continuing to build $G(r/s)$. This simply means
    the rows of $\lambda(r'/s')$ will not be changed. Now that we have added an extra $U$, the next part of $W(r/s)$ is $R^{a_m}$.
    This means we add another row on top of $\lambda(r'/s')$ whish has $a_m$ more boxes than the previous row. By induction,
    the previous row has $\sum_{i=1}^{m-1} a_i$ boxes, and so this new row has $\sum_{i=1}^m a_i$ boxes. Finally, $W(r/s)$ ends
    with $U^{b_m-1}$, and so there are a total of $b_m$ rows at the top, all of this length.

    \bigskip

    $(b)$ The statement says $\mu$ is obtained from $\lambda$ by subtracting one from all parts (i.e. $\mu_i = \lambda_i-1$), and then
    the multiplicities are the same for all but the first (i.e. the greatest) part, which is one less. To see this, consider a part of
    the word $W(r/s)$ of the form $R^{a}U^{b}$. When the snake graph goes right $a$ times, this corresponds to the bottom-most row
    of of the Young diagram of some length $\lambda_i$. After this, every time it goes up (except the last), we get another row of length $\lambda_i$. 
    The exception is the last time it goes up, since if the snake graph continues to go right again after this, the final ``up'' will be part of a longer row.

    \medskip

    At each of these ``up'' steps, the corresponding part of $\mu$ will have length $(\lambda_i-1)$ (including the last ``up''). If this part $R^aU^b$ is not
    the end of the snake graph, then by part $(a)$, $b$ is the number of parts of $\lambda$ of length $\lambda_i$. 
    However, in the special case that $U^{b_m-1}$ is the end of
    the snake graph, we only get $b_m-1$ parts of length $\mu_m$ (but there are $b_m$ parts of length $\lambda_m$).
\end {proof}

\bigskip

\begin {prop} \label{prop:youngs_lattice_gen_fn}
    Let $\lambda = \lambda(r/s)$ and $\mu=\mu(r/s)$. The poset $P(r/s)$ is isomorphic to the interval $[\mu, \lambda]$ in Young's lattice, and
    \[ \mathcal{R}(q) = \frac{1}{q^{|\mu|}} \sum_{\mu \leq \nu \leq \lambda} q^{|\nu|} \]
\end {prop}
\begin {proof}
    There is the natural bijection between $P(r/s)$ and the interval $[\mu,\lambda]$ in Young's lattice, where
    a lattice path corresponds to the Young diagram whose boxes lie above the path. Under this bijection, the rank
    of a path corresponding to a partition $\nu$ is $|\nu|-|\mu|$, which is the number of boxes in the snake graph above the path.
    The formula for $\mathcal{R}(q)$ follows from this bijection together with \textbf{Proposition \ref{prop:rank_function}}.
\end {proof}

\bigskip

\section {Schubert Varieties} \label{sec:Schubert}

\subsection {Definitions}

\bigskip

Let $\Bbb{K}$ be a field, and consider the Grassmannian $\mathrm{Gr}_k(n)$ of $k$-planes in $\Bbb{K}^n$. We will
identify $\mathrm{Gr}_k(n)$ with the set of $k \times n$ matrices of rank $k$, modulo left multiplication by elements of $\mathrm{GL}_k$.
The following definitions and results are all well-known (see \cite{fulton} for a standard reference,
although here we use notational conventions as in \cite{postnikov}).

\bigskip

There is a simple bijection between the interval $[\varnothing, (n-k)^k]$ in Young's lattice (i.e. partitions which fit inside
a $k \times (n-k)$ rectangle) and the set $\binom{[n]}{k}$ of $k$-element subsets of $[n] = \{1,2,\dots,n\}$, given as follows.
Recall from \textbf{Section \ref{sec:wt_fns}} the bijection $Y_{n,k} \to P_{n,k}$ that associates to a Young diagram
$\lambda$ the lattice path from $(0,0)$ to $(n-k,k)$ which is formed by the bottom boundary of $\lambda$.
If the path is traversed backwards, and the steps are labelled $1,2,\dots,n$,
then taking the subset of vertical steps gives an element $I_\lambda \in \binom{[n]}{k}$.

\bigskip

\begin {defn}
    The \emph{open Schubert cell} $\Omega^\circ_\lambda \subseteq \mathrm{Gr}_k(n)$ is the subset whose representatives,
    when written in echelon form, have the identity matrix in the columns indexed by $I_\lambda$. The remaining non-zero
    entries (to the right of the pivots) form the shape of $\lambda$ (but backwards).
\end {defn}

\bigskip

From this definition, it is easy to see that $\dim(\Omega^\circ_\lambda) = |\lambda|$. Once the matrices are in echelon
form, there is a clear parameterization by $\Bbb{K}^{|\lambda|}$.

\bigskip

\begin {ex}
    The partition $\lambda = (3,2,2,1)$ fits inside a $3 \times 4$ rectangle (corresponding to $\mathrm{Gr}_4(7)$).
    The corresponding column set is $I_\lambda = \{1,3,4,6\}$. The open Schubert cell $\Omega^\circ_\lambda$ contains all
    matrices whose echelon form has the shape given in \textbf{Figure \ref{fig:schubert_cell}}.
    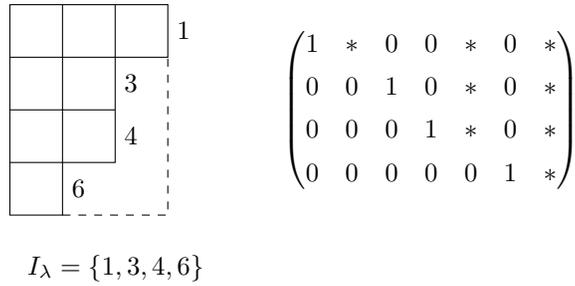
\begin {figure}[h]
    \centering
    \begin {tikzpicture}[scale=0.7]
        \draw (0,0) -- (1,0) -- (1,1) -- (2,1) -- (2,3) -- (3,3) -- (3,4) -- (0,4) -- cycle;
        \draw (1,1) -- (1,4);
        \draw (2,3) -- (2,4);
        \draw (0,1) -- (1,1);
        \draw (0,2) -- (2,2);
        \draw (0,3) -- (2,3);
        \draw[dashed] (1,0) -- (3,0) -- (3,3);

        \draw (3,3.5) node[right] {$1$};
        \draw (2,2.5) node[right] {$3$};
        \draw (2,1.5) node[right] {$4$};
        \draw (1,0.5) node[right] {$6$};

        \draw (2,-1) node {$I_\lambda = \{1,3,4,6\}$};

        \draw (8,2) node {
            $\begin{pmatrix} 1 & \ast & 0 & 0 & \ast & 0 & \ast \\[1ex]
                             0 & 0    & 1 & 0 & \ast & 0 & \ast \\[1ex]
                             0 & 0    & 0 & 1 & \ast & 0 & \ast \\[1ex]
                             0 & 0    & 0 & 0 & 0    & 1 & \ast \end{pmatrix}$
        };
    \end {tikzpicture}
    \caption {(left) The subset $I_\lambda \in \binom{[n]}{k}$ corresponding to a partition $\lambda$, and (right) the echelon form of a matrix representative in $\Omega^\circ_\lambda$}
    \label {fig:schubert_cell}
    \end {figure}
\end {ex}

\bigskip

\begin {defn}
    The \emph{Schubert variety} (or \emph{closed Schubert cell}) is the closure $\Omega_\lambda := \overline{\Omega^\circ_\lambda}$.
\end {defn}

\bigskip

It is a well-known fact that $\Omega_\lambda$ is the disjoint union of all open Schubert cells corresponding to partitions that
fit inside $\lambda$:
\[ \Omega_\lambda = \bigcup_{\mu \leq \lambda} \Omega^\circ_\mu \]
In particular, the Grassmannian is the closure of the ``biggest'' cell $\Omega^\circ_\lambda$, where $\lambda = (n-k)^k$ is the entire $k \times (n-k)$ 
rectangle, and thus $\mathrm{Gr}_k(n)$ is the disjoint union of all the open Schubert cells. 

\bigskip

Now, putting this all together, we are in a position to see the connection between the geometric and combinatorial
interpretations of the $q$-binomial coefficients given in \textbf{Theorem \ref{thm:finite_variety_sizes}(b)} and \textbf{Theorem \ref{thm:q_binom_gen_fn}}.
If we take $\Bbb{K} = \Bbb{F}_q$, then clearly $\left| \Omega^\circ_\lambda \right| = q^{|\lambda|}$. The fact that the Grassmannian is
the disjoint union of all open Schubert cells implies that
\[ \left| \mathrm{Gr}_k(n) \right| = \sum_{\lambda \leq (n-k)^k} \left| \Omega^\circ_\lambda \right| = \sum_{\lambda \leq (n-k)^k} q^{|\lambda|} = \binom{n}{k}_q \]

\bigskip

\subsection {Numerators of $q$-Rationals}

\bigskip

Finally, we return to the discussion of $\mathcal{R}(q)$, then numerator of the $q$-rational $\left[ \frac{r}{s} \right]_q$, and
we will show that it is the size of a certain subvariety of $\mathrm{Gr}_k(n)$ (over $\Bbb{F}_q$),
where $k$ and $n$ can be determined by the continued fraction $\frac{r}{s} = [a_1,a_2,\dots,a_{2m}]$.
Let $n = \sum_{i=1}^{2m} a_i$ and $k = \sum_{i=1}^m a_{2i}$.

\bigskip

\begin {thm} \label{thm:main}
    Let $\frac{r}{s} = [a_1,\dots,a_{2m}] \in \Bbb{Q}$, and let $\left[\frac{r}{s}\right]_q = \frac{\mathcal{R}(q)}{\mathcal{S}(q)}$ be the corresponding $q$-rational.
    Let $\mu = \mu(r/s)$ and $\lambda = \lambda(r/s)$. 
    Then up to a factor of $q^{|\mu|}$, the polynomial $\mathcal{R}(q)$ is the number of $\Bbb{F}_q$ points in the union of open Schubert
    cells in the interval $[\mu,\lambda]$ in the Grassmannian $\mathrm{Gr}_k(n)$:
    \[ q^{|\mu|} \cdot \mathcal{R}(q) = \left| \bigcup_{\mu \leq \nu \leq \lambda} \Omega^\circ_\nu \right| \]  
\end {thm}
\begin {proof}
    Recall from \textbf{Proposition \ref{prop:youngs_lattice_gen_fn}} that the numerator $\mathcal{R}(q)$ is the
    rank generating function of the interval $[\mu,\lambda]$ in Young's lattice (scaled by $q^{|\mu|}$).
    The width of the snake graph (also the width of $\lambda$) is $\sum_{i=1}^m a_{2i-1}$, and the height is $\sum_{i=1}^m a_{2i}$.
    Since the parameters $n$ and $k$ of the Grassmannian correspond to Young diagrams in a $k \times (n-k)$ rectangle,
    this implies that $k = \sum_{i=1}^m a_{2i}$ and $n-k = \sum_{i=1}^m a_{2i-1}$.

    \medskip

    Since $q^{|\nu|}$ is the size of the open Schubert cell $\Omega^\circ_\nu$, we see that
    \[ \mathcal{R}(q) = \frac{1}{q^{|\mu|}} \sum_{\mu \leq \nu \leq \lambda} \left| \Omega^\circ_\nu \right| \]
\end {proof}

\bigskip

\begin {ex}
    If $\frac{r}{s} = \frac{n}{1}$ or $\frac{r}{s} = \frac{n+1}{n}$, then the snake graph $G(r/s)$ is a straight row or column of boxes.
    In this case, the partition $\mu$ is empty, and $\lambda$ is the full $n \times 1$ or $1 \times n$ rectangle. In this case \textbf{Theorem \ref{thm:main}}
    reduces to a special case of \textbf{Theorem \ref{thm:finite_variety_sizes}(b)}, which, as described in \textbf{Example \ref{ex:projective_space}},
    simply says that $|\Bbb{P}^{n-1}| = [n]_q$. The interpretation given in \textbf{Theorem \ref{thm:main}} in terms of Schubert cells just states
    that $\Bbb{P}^{n-1} = \bigcup_{k=1}^{n-1} X_k$, where $X_k$ consists of those points whose homogeneous coordinates have the form
    \[ [x_0:x_1: \cdots : x_{n-1}] = [0:0: \cdots : 0 : 1 : \ast : \cdots : \ast] \]
    where $x_k=1$ and $x_i = 0$ for $i < k$.
\end {ex}

\bigskip

\begin {ex}
    Continuing with \textbf{Example \ref{ex:seven_thirds}} and \textbf{Example \ref{ex:seven_thirds_again}}, let $\frac{r}{s} = \frac{7}{3} = [2,3]$. 
    The corresponding $q$-rational is
    \[ \left[ \frac{7}{3} \right]_q = \frac{\mathcal{R}(q)}{\mathcal{S}(q)} = \frac{1+2q+2q^2+q^3+q^4}{1+q+q^2} \]
    The snake graph $G(7/3)$ has word $W(G) = RUU$, and is pictured (along with the poset $P(7/3)$) in \textbf{Figure \ref{fig:seven_thirds_again}}. The snake graph
    has width 2 and height 3, which means the corresponding union of Schubert cells lives in $\mathrm{Gr}_3(5)$ (since $k=3$ and $n-k=2$).

    \medskip

    In the notation of \textbf{Theorem \ref{thm:main}}, we have $\lambda = (2,2,2)$ and $\mu = (1,1)$.
    In  this case \textbf{Theorem \ref{thm:main}} says that $q^2 \mathcal{R}(q)$ counts the number of points in $\bigcup_{\nu \in [\mu,\lambda]} \Omega^\circ_\nu$,
    which consists of all 3-dimensional subspaces of $\Bbb{F}_q^5$ which have a matrix representative of one of the following seven forms:
    \[
        \begin{pmatrix}
            1 & 0 & 0 & \ast & \ast \\[1ex]
            0 & 1 & 0 & \ast & \ast \\[1ex]
            0 & 0 & 1 & \ast & \ast 
        \end{pmatrix}, \quad
        \begin {pmatrix}
            1 & 0 & \ast & 0 & \ast \\[1ex]
            0 & 1 & \ast & 0 & \ast \\[1ex]
            0 & 0 & 0    & 1 & \ast
        \end{pmatrix}, \quad
        \begin{pmatrix}
            1 & 0 & \ast & \ast & 0 \\[1ex]
            0 & 1 & \ast & \ast & 0 \\[1ex]
            0 & 0 & 0    & 0    & 1
        \end{pmatrix}, \quad
        \begin{pmatrix}
            1 & \ast & 0 & 0 & \ast \\[1ex]
            0 & 0    & 1 & 0 & \ast \\[1ex]
            0 & 0    & 0 & 1 & \ast
        \end{pmatrix},
    \]
    \[
        \begin{pmatrix}
            0 & 1 & 0 & 0 & \ast \\[1ex]
            0 & 0 & 1 & 0 & \ast \\[1ex]
            0 & 0 & 0 & 1 & \ast 
        \end{pmatrix}, \quad
        \begin{pmatrix}
            1 & \ast & 0 & \ast & 0 \\[1ex]
            0 & 0    & 1 & \ast & 0 \\[1ex]
            0 & 0    & 0 & 0    & 1 
        \end{pmatrix}, \quad
        \begin{pmatrix}
            0 & 1 & 0 & \ast & 0 \\[1ex]
            0 & 0 & 1 & \ast & 0 \\[1ex]
            0 & 0 & 0 & 0    & 1 
        \end{pmatrix}
    \]
\end {ex}

\bigskip

\section* {Acknowledgments}

\bigskip

I would like to acknowledge the support of the NSF grant DMS-1745638. I would also like to thank
Valentin Ovsienko and Vic Reiner for helpful comments and conversations.

\bigskip

\bibliographystyle{alpha}
\bibliography{qrats}

\end {document}